\documentclass[a4paper,12pt]{article}
\usepackage{inputenc}
\usepackage[english]{babel}
\usepackage{amsfonts}
\usepackage{graphics}
\usepackage{wrapfig}
\usepackage[pdflatex]{graphicx}
\DeclareGraphicsExtensions{.pdf,.png}
\date{}
\newtheorem{Th}{Theorem}[section]

\begin{document}
\renewcommand{\baselinestretch}{0.96}
\renewcommand{\thesection}{\arabic{section}}
\renewcommand{\theequation}{\thesection.\arabic{equation}}
\csname @addtoreset\endcsname{equation}{section}
 \sloppy
\vspace*{1cm}

\large

\begin{center} {\bf Numerical solution of linear differential equations
with discontinuous coefficients and Henstock integral}\\
S.F. Lukomskii, D.S. Lukomskii
\end{center}
\begin{abstract}In this article we consider the problem  of approximative solution  of linear differential equations $y'+p(x)y=q(x)$
with discontinuous coefficients $p$ and $q$. We assume that  coefficients of  such equation are Henstock integrable functions. 
To find  the approximative solution we change the original Cauchy  problem to another problem with piecewise-constant coefficients. 
The sharp solution of this new problems is the approximative solution of the original Cauchy  problem. 
We  find the degree  approximation in terms of modulus of continuity $\omega_\delta (P),\ \omega_\delta (Q)$, 
where $P$ and $Q$ are $f$-primitive for coefficients $p$ and $q$.

 \end{abstract}

\section{Introduction}
In the classical initial value problem for a linear differential equation of the first order
\begin{equation} \label{eq1.1}
y\prime+p(x)y=q(x),\quad y(a)=y_0,\quad x\in [a,b],
\end{equation}
 the coefficients $p(x)$ and $q(x)$ are continuous functions.
 However, some problems of dry friction  and electric circuit with relay are given the equations with the discontinuous functions 
$p$ and $q$. For example the  RL-electric circuit with relay is described by a linear differential equation
 $$
 \frac{di}{dt}+\frac{R(t)}{L}i=\frac{e(t)}{L}
 $$
   with the discontinuous function $R(t)$.
 In this case, it is assumed that the functions $p(x)$ and $p(x)$ are (L)-integrable and a function $y(x)$ is called a solution to
 equation (\ref{eq1.1}) if $y(x)$ is absolutely continuous and satisfies the equation (\ref{eq1.1}) almost everywhere on $[a,b]$.

 There are no effective methods for the approximate solution of equations with unbounded coefficients $p(x)$ and $q(x)$. 
If the coefficients $p(x)$ and $q(x)$ are unbounded  in some neighborhood of the point $a$ then Ruge-Kutta method does not work.
  If the coefficients $p(x)$ and $q(x)$ are unbounded  in some neighborhood of the interior point $c\in (a,b)$ then
 Runge-Kutta method  has a very large error, usually more than 1.

 Some authors use  Haar and Walsh functions to solve linear equations \cite{OhkKob},\cite{RazOrd1},\cite{RazOrd2}. 
In \cite{GGRT1}, \cite{GGRT2} G. Gat and R. Toledo propose to approach the  solution $y(x)$ by the Walsh polynomial
 $$\tilde y_n(x)=\sum_{k=0}^{2^n-1}c_kW_k(x).$$ In \cite{GGRT1} for continue $q(x), (x\in[0,1])$ and 
 $p(x)=const$   an estimate for the error $|y(x)-\tilde y(x)|$ is obtained.
  In \cite{GGRT2} the authors consider the case when  $q\in L(0,1)$ is a continuous function
 on $[0,1[$ and prove that $\tilde y_n(x)$ converges uniformly to the solution $y(x)$ on the interval $[0,1[$.

  In \cite{LLT}, the authors present the derivative $y'$ of the solution $y$ as a Haar expansion and obtain an estimate
 of the approximate solution in terms of the modulus of continuity of the coefficients $p(x)$ and $q(x)$.
 This method can also be used for equations with unbounded coefficients $p(x)$ and $q(x)$.

In this article we will assume that $p(x)$ and $q(x)$ are Henstock integrable  functions on the interval $[a,b]$.
 We construct the  approximate solution $\tilde y(x)$ and obtain the estimate of the error $|y(x)-\tilde y(x)|$
 in terms of modulus of continuity $\omega_{\frac{1}{2^n}}(e^{P})$, $\omega_{\frac{1}{2^n}}(e^{-P})$,  
and $\omega_{\frac {1}{2^n}}(Q)$, where $P$ and $Q$ are $f$-primitive for $p$ and $q$ respectively.

The paper is organized as follows.
In Sec. 2, we recall some facts from Henstock integral.
In Sec. 3, we indicate the necessary  and sufficient condition, under which  the Cauchy problem has a solution.
This solution is given in terms of the Hanstock integral.
In Sec. 4, we construct the approximative solution  and find the  error.
 In Sec. 5, we give two examples.

\section{Henstock integral on the interval}

Any function  $\delta(x)>0$ on $[a,b]$ is said to be a gauge.
 Let $\mathfrak X=(x_k)_{k=0}^n$ be a partition of the interval $[a,b]$. The  point $\xi_k\in[x_{k-1},x_k]$  is called a tag of $[x_{k-1},x_k]$ , the set of ordered pairs   $([x_{k-1},x_k],\xi_k)_{k=1}^n$ is called a tagged partition and is denote by
$ \stackrel{\circ}{\mathfrak X}= ([x_{k-1},x_k],\xi_k)_{k=1}^n$.

The tagged partition $ \stackrel{\circ}{\mathfrak X}=
([x_{k-1},x_k],\xi_k)_{k=1}^n$ of the interval
$[a,b]$  is called  $\delta$-fine and is denote by $\stackrel{\circ}{\mathfrak X}\ll\delta$ if for any
$k=1,\dots,n$
$$
  |x_{k-1}-x_k|<\delta(\xi_k).
$$
  It is known that for any gauge  $\delta(x)>0$ on $[a,b]$  there exists a $\delta$-fine partition $ \stackrel{\circ}{\mathfrak X}=
([x_{k-1},x_k],\xi_k)_{k=1}^n$ of $[a,b]$.

 A function $f:[a,b]\rightarrow \mathbb R$
is said to be  Henstock-integrable ( or generalized Riemann-integrable)  on the interval $[a,b]$, 
if there exists a number  $I(f)\in\mathbb R$ such that
$$
\forall\;\varepsilon>0\;\exists\;\delta(x)>0\;\mbox{on}\;[a,b]\;\forall\;
\stackrel{\circ}{\mathfrak
X}\ll\delta(x),\quad|S(\stackrel{\circ}{\mathfrak
 X},f)-I(f)|<\varepsilon.
$$
The number  $I(f)$ is called Henstock integral and  is denoted by
$(R^*)\int\limits_a^b f(x)\,dx$ or  $\int\limits_a^b f(x)\,dx$.

 The collection of all functions that are Henstock integrable on $[a,b]$ is denoted by $R^*(a,b)$.
A function  $f:[a,b]\rightarrow \mathbb R$ is called absolutely integrable if
 $f\in R^*(a,b)$ and $|f|\in R^*(a,b)$. There exists Henstock integrable functions that are not absolutely integrable.
 If the function $f:[a,b]\rightarrow \mathbb R$ is absolutely integrable then $f\in L(a,b)$

 The function $G:[a,b]\rightarrow R$ is called a  $c$-primitive ($f$-primitive) for a function  $g$ if 
 $G$ is continuous on $[a,b]$ and there exist a countable (finite) set $E\subset  [a,b]$ such
that  $G'(x)=g(x)$ on $[a,b]\backslash E$ . We will use next properties of Henstock integral.
\begin{Th}[\cite{BG}]
 If $f:[a,b]\rightarrow R $ has a $c$-primitive $F$ with a exceptional set $E$ , then  $f\in R^*(a,b)$ and for all $x$
 $$
 \int_a^xf(t)dt=F(x)-F(a).
 $$
\end{Th}
It follows that for  $x\in [a,b]\setminus E$
$$
 \frac{d}{dx}\int_a^xf(t)dt=f(x).
 $$
 \begin{Th}[\cite{BG}]
    Let $f\in R^*(a,b)$ and  $F(x)=\int_a^xf(t)dt$. The function  $f$ is absolutely integrable on $ [a,b]$ if and only if 
$\bigvee_a^b(F)<+\infty$. In this case,
 $$
 \int_a^b|f(t)|dt=\bigvee_a^b(F).
 $$
\end{Th}

\begin{Th}[\cite{BG}]
 If  $f\in R^*(a,b)$ and  $g$ is monotone on $[a,b]$, then there exists $\xi\in[a,b]$ such that
 $$
 \int_a^bf(x)g(x)dx=g(a)\int_a^{\xi}f(x)dx+g(b)\int_{\xi}^{b}f(x)dx.
 $$
\end{Th}
\begin{Th}[\cite{BG}]
    Let $F$ and $G$ be a c-primitives  on $[a,b]$. Then $F'G\in R^*(a,b)$ if and only if $FG'\in R^*(a,b)$. In this case,
 $$
 \int_a^bF'Gdt= \Biggl. F(t)G(t)\Biggr|_a^b -\int_a^bFG'dt.
 $$
\end{Th}
\begin{Th}[\cite{BG}, Hake`s theorem]
Let $f:[a,b]\rightarrow \mathbb R$ and $f\in R^*(a,c)$ for any $c\in (a,b)$. Then $f\in R^*(a,b)$ if and only if there exists
$$
\lim_{c\rightarrow b-0}\int_a^cf(x)dx =I.
$$
In this case, $I=(R^*)\int_a^bf(x)dx$.
\end{Th}

 A detailed exposition of the Henstock integral theory can be found in   \cite{BG},\cite{RA}.

\section{Linear differential equations and Henstock integral}
   Let $p,q:[a,b]\rightarrow\mathbb R$ be two  continuous functions that
     differentiable on the interval $[a,b]$, with the exception of a countable set
   $E$.
 We will consider   the classical Cauchy initial value problem
\begin{equation} \label{eq3.1}
y'+p'(x)y=q'(x),\quad x\in[a,b]\setminus E ,
\end{equation}
\begin{equation} \label{eq3.2}
y(a)=y_0.
\end{equation}
 Theorem 2.1 follows that functions  $p'(x)$ and $q'(x)$ are Henstock integrable.
 This is a weaker condition than $p',q'\in L(a,b)$.\\
{\bf Example 1.} Define the function $q$ on $[a,b]$ in the following way. Let $x_n=a+\frac{b-a}{2^n}$. 
Assume $q(a)=q(x_n)=0, q(\frac{x_n+x_{n+1}}{2})=\frac{1}{n}$, $q(x)$ 
is lineal on $[x_{n+1},\frac{x_n+x_{n+1}}{2}]$ and $[\frac{x_n+x_{n+1}}{2},x_n]$. 
Then $q'\in R^*(a,b)$ but $q'\notin L(a,b)$.

\begin{Th} {\it Equation (\ref{eq3.1}) has a continuous solution  that is differentiable on the set  $[a,b]\setminus E$
if and only if  the function  $e^{p(x)}q'(x)$ has a $c$-primitive differentiable on $[a,b]\setminus E$.}
\end{Th}
{\bf Necessity.} Let $y(x)$ be a solution of (\ref{eq3.1}), that is
$$
y'(x)+p'(x)y(x)=q'(x)
$$
for all $x\in [a,b]\setminus E$. Then
$$
e^{p(x)}q'(x)+p'(x)y(x)e^{p(x)}=q'(x)e^{p(x)}
$$
for all  $x\in [a,b]\setminus E$
or in another words
\begin{equation} \label{eq3.3}
(y(x)e^{p(x)})'=q'(x)e^{p(x)}\qquad (x\in [a,b]\setminus E)
\end{equation}
It means, that the function  $q'(x)e^{p(x)}$ has $c$-primitive  $y(x)e^{p(x)}$.\\
{\bf Sufficiently.} Let $q'(x)e^{p(x)}$ has  $c$-primitibe $F(x)$, that is
$$
F'(x)=q'(x)e^{p(x)}\qquad x\in[a,b]\setminus E.
$$
 Let us denote $y(x)=\frac{F(x)}{e^{p(x)}}\Leftrightarrow F(x)=y(x)e^{p(x)}$  $(x\in [a,b]\setminus E)$. Then
$$
\forall x\in [a,b]\setminus E \quad y'(x)e^{p(x)}+y(x)e^{p(x)}p'(x)=q'(x)e^{p(x)}\Leftrightarrow
$$
$$
y'(x)+y(x)p'(x)=q'(x). \square
$$
{\bf Corollary.} {\it A solution  of Cauchy initial value problem (\ref{eq3.1})-(\ref{eq3.2}) is given by the formula

$$
y(x)=e^{p(a)-p(x)}g(a)+e^{-p(x)}\int_a^xq'(t)e^{p(t)}d\,t ,
$$
where the integral is an Henstock integral.}\\
{\bf Proof.} Equality (\ref{eq3.3}) follows, that the function  $y(x)e^{p(x)}$ is $c$-primitive for  $q'e^{p(x)}$, 
it means  $q'e^{p(x)}$ is Henstock integrable end the equality
$$
\int_a^xq'(t)e^{p(t)}d\,t=y(x)e^{p(x)}-y(a)e^{p(a)}.
$$
holds. $\qquad \square$

{\bf Example 2.} It is possible to construct the continuous functions $p$ and $q$  so that the function
 $q'(x)e^{p(x)}$ has a c-primitive, but $q'(x)e^{p(x)}\notin L(a,b)$.
For simplicity, we consider the case $[a,b]=[0,1]$  and select the function $q(x)$ as in Example 1. 
In this case $x_n=2^{-n}, \ q(x_n)=q(0)=0$,  $q(x)$ is lineal on $[x_{n+1},\xi_n]$ and $[\xi_n,x_n]$,
 where $\xi_n=\frac12(x_x+x_{n+1})$.  Now we define the function $p(x)$ from the conditions:\\
(a)$e^{p(2^{-n})}=\beta_n>1, \beta_n\downarrow 1 \ (n\to \infty)$,\\
(b)$e^{p(x)}$ is lineal on any interval $[2^{-n-1},2^{-n}]$.
It is evident that the series
$$
\sum_{n=1}^\infty \int_{2^{-k+1}}^{2^{-k}}q'(x)e^{p(x)}dx
$$
converges. It follows from the Hake theorem that    $f(x)=q'(x)e^{p(x)}\in R^*(0,1)$. Therefore 
the function $F(x)=\int_0^xf(t)dt$ is continuous. Since the function $f(x)=q'(x)e^{p(x)}$ 
is continuous on   any interval $(2^{-n-1},2^{-n})$, it follows that $F'(x) =q'(x)e^{p(x)}$ on 
  any interval $(2^{-n-1},2^{-n})$. It means that $F(x)$ is c-primitive for $q'(x)e^{p(x)}$.
It is not difficult to check that  $f(x)=q'(x)e^{p(x)}\notin L(0,1)$.

\section{Approximate solution of Cauchy  problem (\ref{eq3.1})-(\ref{eq3.1}) on interval [0,1]}
Now we will find an approximate solution of Cauchy initial value problem
\begin{equation} \label{eq4.1}
y'+p'(x)y=q'(x),\quad x\in[0,1]\setminus E ,
\end{equation}
\begin{equation} \label{eq4.2}
y(0)=y_0.
\end{equation}
  We assume that the functions $p$ and $q$ are continuous and have derivatives
  with the exception of
  some countable set $E$. We also assume, that $e^{p(x)}q'(x)$ has a $c$-primitive differentiable on $[a,b]\setminus E$.

Choose an arbitrary $n\in \mathbb N$, define the functions  $\tilde p(x)$ and $\tilde q(x)$ by equalities
$$
\tilde p\Biggl(\frac{k}{2^n}\Biggr)=p\Biggl(\frac{k}{2^n}\Biggr), \tilde q\Biggl(\frac{k}{2^n}\Biggr)=
q\Biggl(\frac{k}{2^n}\Biggr),
$$
$$
\tilde p(x)=p\Biggl(\frac{k}{2^n}\Biggr)+2^n\Biggl(x-\frac{k}{2^n}\Biggr)\Biggl(p\Biggl(\frac{k+1}{2^n}\Biggr)
-p\Biggl(\frac{k}{2^n}\Biggr)\Biggr), x\in \Biggl[\frac{k}{2^n},\frac{k+1}{2^n}\Biggr],
$$
$$
\tilde q(x)=q\Biggl(\frac{k}{2^n}\Biggr)+2^n\Biggl(x-\frac{k}{2^n}\Biggr)\Biggl(q\Biggl(\frac{k+1}{2^n}\Biggr)
q\Biggl(\frac{k}{2^n}\Biggr)\Biggr), x\in \Biggl[\frac{k}{2^n},\frac{k+1}{2^n}\Biggr],
$$
and consider the Cauchy initial value problem
\begin{equation} \label{eq4.3}
\tilde y'+\tilde p'\tilde y=\tilde q'\ ,
\end{equation}
\begin{equation} \label{eq4.4}
\tilde y(0)=y_0.
\end{equation}
It is evident that the function  $e^{\tilde p(x)}\tilde q'$ has a  $f$-primitive.
By theorem 3.1 the functions
$$
y(x)=y_0e^{p(0)-p(x)}+e^{-p(x)}\int\limits_0^xq'(t)e^{p(t)}d\,t,
$$
$$
\tilde y(x)=y_0e^{\tilde p(0)-\tilde p(x)}+e^{-\tilde p(x)}\int\limits_0^x\tilde q'(t)e^{\tilde p(t)}d\,t.
$$
are solutions of  Cauchy problems (\ref{eq4.1})-(\ref{eq4.2}) and (\ref{eq4.3}) - (\ref{eq4.4}) respectively.
 The function $\tilde y(x)$ is the approximate solution of  Cauchy problem (\ref{eq4.1})-(\ref{eq4.2}). 
In the following theorem, we indicate an estimate for  the distance
 $y(x)-\tilde y(x)$.\\
\begin{Th} The following inequality
$$
|y(x)-\tilde y(x)|\le
C_{-1}\omega_{\frac{1}{2^n}}(e^{-p})+C_{1}\omega_{\frac{1}{2^n}}(e^{p})+\omega_{\frac{1}{2^n}}(q)C_0+C_2 \omega_{\frac{1}{2^n}}(q)\omega_{\frac{1}{2^n}}(p).
$$
holds, where
$$
C_{-1}=|y_0|e^{p(0)}+\|e^p\|_{C[0,1]}\bigvee\limits_0^1 q,\quad
C_{1}=2\|e^p\|_{C[0,1]}\bigvee\limits_0^1 q,
$$
$$
C_0=\|e^p\|_{C[0,1]}+\|e^p\|_{C[0,1]}\bigvee\limits_0^1 e^p,\quad C_2=\|e^p\|_{C[0,1]}^2.
$$
\end{Th}
{\bf Proof}
1)First we estimate the difference  $y(x)-\tilde y(x)$ for $x=\frac{k}{2^n}, k=0,1,...,2^n$. We have
$$
y(x)-\tilde y(x)=e^{-p(\frac{k}{2^n})}
\int\limits_0^{\frac{k}{2^n}}(q'(t)e^{p(t)}-\tilde q'(t)e^{\tilde p(t)})d\,t=
$$
$$
e^{-p(\frac{k}{2^n})}\int_0^{\frac{k}{2^n}}(q'(t)-\tilde q'(t))e^{ p(t)}d\,t+
e^{-p(\frac{k}{2^n})}\int\limits_0^{\frac{k}{2^n}} \tilde q'(t)(e^{p(t)}-e^{\tilde p(t)})d\,t=I_1+I_2.
$$

     To estimate integrals in $I_1$ and  $I_2$ we will assume that  $p'$ and  $q'$ -- are Henstock absolutely integrable.\\
  Assume $I_1$. Integrating by parts and using  the equality $q'(\frac{j}{2^n})=\tilde q'(\frac{j}{2^n})$ we obtain
$$
\left|\int\limits_0^{\frac{k}{2^n}}(q'(t)-\tilde q'(t))e^{p(t)}d\,t\right|\le \left|q(t)-\tilde q(t)|_0^\frac{k}{2^n}\right|+ \left|\int\limits_0^{\frac{k}{2^n}}(q(t)-\tilde q(t))(e^{p(t)})'d\,t\right|\le
$$
$$
\le\omega_{\frac{1}{2^n}}(q)\int\limits_0^{\frac{k}{2^n}}|(e^{p(t)})'|d\,t\le
\omega_{\frac{1}{2^n}}(q)\bigvee\limits_0^1e^{p(\cdot)}.
$$
So
$$
|I_1|\le \|e^{-p(\cdot)}\|_{C(0,1)}\omega_{\frac{1}{2^n}}(q(\cdot))\bigvee\limits_0^1e^{p(\cdot)}.
$$
   Since the function $e^{\tilde p(t)}$ is monotonic on any interval  $[\frac{j}{2^n},\frac{j+1}{2^n}]$, it follow that
$|(e^{p(t)}-e^{\tilde p(t)})|\le \omega_{\frac{1}{2^n}}(e^{p(\cdot)})$. Therefore
$$
|I_2|\le \|e^{-p(\cdot)}\|\cdot \omega_{\frac{1}{2^n}}(e^{p(\cdot)})\bigvee\limits_0^1q(\cdot),
$$
and
$$
\left| y\left(\frac{k}{2^n}\right)-\tilde y\left(\frac{k}{2^n}\right)\right|\le 
\|e^{-p(\cdot)}\|_{C[0,1]}\left( \omega_{\frac{1}{2^n}}(q(\cdot)) \bigvee\limits_0^1e^{p(\cdot)}+\omega_{\frac{1}{2^n}}(e^{p(\cdot)})\cdot \bigvee\limits_0^1q(\cdot) \right).
$$
2) Now we consider the case  $x\in \left[\frac{k}{2^n},\frac{k+1}{2^n}\right]$. Let us write the difference 
$y(x)-\tilde y(x)$ in the form
$$
y(x)-\tilde y(x)=y_0 e^{p(0)}(e^{-p(x)}-e^{-\tilde p(x)})+
$$
$$
+(e^{-p(x)}-e^{-\tilde p(x)})\left(\int\limits_0^{\frac{k}{2^n}}q'(t)e^{p(t)}d\,t+
\int\limits_{\frac{k}{2^n}}^xq'(t)e^{p(t)}d\,t\right)+
$$
$$
+e^{-\tilde p(x)}\left(\int\limits_0^{\frac{k}{2^n}}q'(t)e^{p(t)}d\,t-
\int\limits_0^{\frac{k}{2^n}}\tilde q'(t)e^{\tilde p(t)}d\,t \right)+
$$
\begin{equation} \label{eq4.5}
+e^{-\tilde p(x)}\left(\int\limits^{\frac{k}{2^n}}_xq'(t)e^{p(t)}d\,t-
\int\limits_{\frac{k}{2^n}}^x\tilde q'(t)e^{\tilde p(t)}d\,t \right)=
A_1+A_2+(A_3+A_4)e^{-\tilde p(x)}.  
                                             \end{equation}
We will estimate $A_l$ $(l=1,2,3,4)$.\\
1) Since the function  $e^{-\tilde p(x)}$ is monotonic on any interval  $[\frac{j}{2^n},\frac{j+1}{2^n}]$, it follow that
\begin{equation} \label{eq4.6}
|e^{-p(x)}-e^{-\tilde p(x)}|\le\omega_{\frac{1}{2^n}}(e^{-p}).         
               \end{equation}
2) Using again (\ref{eq4.6}), we get
\begin{equation} \label{eq4.7}
|A_2| \le \omega_{\frac{1}{2^n}}(e^{-p})\left|\int\limits_0^x q'(t)e^{p(t)}d\,t\right|\le
\omega_{\frac{1}{2^n}}(e^{-p})\|e^p\|_{C[0,1]}\bigvee\limits_0^1q.                    \end{equation}
3) An estimate for $A_3$ was obtained earlier
\begin{equation} \label{eq4.8}
|A_3|\le \left((\omega_{\frac{1}{2^n}}(q))\bigvee\limits_0^1 e^p+\omega_{\frac{1}{2^n}}(e^p)\bigvee\limits_0^1 q\right).
\end{equation}
4)Let us write  $A_4$ in the form
\begin{equation} \label{eq4.9}
A_4=\int\limits_{\frac{k}{2^n}}^xq'(t)(e^{p(t)}-e^{\tilde p(t)})d\,t+\int\limits_{\frac{k}{2^n}}^xe^{\tilde p(t)}(q'(t)\tilde q')d\,t
\end{equation}
Since the function $e^{\tilde p(t)}$ is monotonic on the interval $[\frac{k}{2^n},\frac{k+1}{2^n}]$, both integrals exist.
 For the first integral, we have the obvious inequality
$$
\int\limits_{\frac{k}{2^n}}^xq'(t)(e^{p(t)}-e^{\tilde p(t)})d\,t\le 
\omega_{\frac{1}{2^n}}(e^p)\cdot \int\limits_{\frac{k+1}{2^n}}^{\frac{k}{2^n}}|q'(t)|d\,t\le \omega_{\frac{1}{2^n}}(e^p)\bigvee\limits_0^1q.
$$
Integrating by parts in the second integral in (\ref{eq4.9}) we have
$$
\left|\int\limits_{\frac{k}{2^n}}^x e^{\tilde p(t)}(q'(t)-\tilde q'(t))d\,t\right|\le |q(x)-\tilde q(x)|+\omega_{\frac{1}{2^n}}(q) \int\limits_{\frac{k}{2^n}}^{\frac{k+1}{2^n}}|(e^{\tilde p(t)})|d\,t\le
$$
$$
\le \omega_{\frac{1}{2^n}}(q)+\omega_{\frac{1}{2^n}}(q)\left|
\int\limits_{\frac{k}{2^n}}^{\frac{k+1}{2^n}}e^{\tilde p(t)}\cdot \tilde p'(t)d\,t \right|\le
$$
$$
\le \omega_{\frac{1}{2^n}}(q)+\omega_{\frac{1}{2^n}}(q)\|e^p\|_{C[0,1]}\left|p\left(\frac{k+1}{2^n}\right)-p\left(\frac{k}{2^n}\right)\right|\le
$$
$$
\le \omega_{\frac{1}{2^n}}(q)(1+\|e^p\|_{C[0,1]}\omega_{\frac{1}{2^n}}(p)).
$$
Finally, we obtain
\begin{equation} \label{eq4.10}
|A_4|\le \omega_{\frac{1}{2^n}}(e^p)\bigvee\limits_0^1 q+\omega_{\frac{1}{2^n}}(q)(1+\|e^p\|_{C[0,1]}\omega_{\frac{1}{2^n}}(p)).
\end{equation}
Substituting  inequalities   (\ref{eq4.6})- (\ref{eq4.8})  and (\ref{eq4.10})  in (\ref{eq4.5}) , we get
$$
|y(x)-\tilde y(x)|\le |y_0|e^{p(0)}\omega_{\frac{1}{2^n}}(e^{-p})+\omega_{\frac{1}{2^n}}(e^{-p})\|e^p\|_{C[0,1]}\bigvee\limits_0^1 q+
$$
$$
+\|e^p\|_{C[0,1]}\left(\omega_{\frac{1}{2^n}}(q)\bigvee\limits_0^1 e^p+
\omega_{\frac{1}{2^n}}(e^{p})\bigvee\limits_0^1 q\right)+
$$
$$
+\|e^p\|_{C[0,1]}\left(\omega_{\frac{1}{2^n}}(e^p)\bigvee\limits_0^1 q+
\omega_{\frac{1}{2^n}}(q)(1+\|e^p\|_{C[0,1]} \omega_{\frac{1}{2^n}}(p))\right)=
$$
$$
=C_{-1}\omega_{\frac{1}{2^n}}(e^{-p})+C_{1}\omega_{\frac{1}{2^n}}(e^{p})+\omega_{\frac{1}{2^n}}(q)C_0+C_2 \omega_{\frac{1}{2^n}}(q)\omega_{\frac{1}{2^n}}(p).\quad\square
$$

\section{Some examples}
{\bf Example 3.} Let us consider the Cauchy problem
\begin{equation} \label{eq5.1}
\left\{
\begin{array}{ll}
y'+\frac{1}{2\sqrt{x}}y=1+\frac{1}{\sqrt{x}},&\quad x\in [0,1]\\
y(0)=0.
\end{array}
\right.
\end{equation}
Here $p(x)=\sqrt{x},\ q(x)=x+2\sqrt{x}$. The solution ${y(x)=2\sqrt{x}}$ of this problem is a continuous function on $[0,1]$,
 but the derivative $y'(0)$ not exists.  Denote by  $\tilde y(x)$ the approximative solution for some $N=2^n>1$.
 In the table 1 we give the approximative solution of Cauchy problem (\ref{eq5.1}).
\vspace{0.5cm}

\hspace{0.5cm} \begin{tabular}{|l|l|l|l|l|l|}\hline
 &N=&16&32&64&128\\ \hline
x&y(x)&$\tilde y(x)$&$\tilde y(x)$&$\tilde y(x)$&$\tilde y(x)$\\ \hline
0& 0&0 &0&0&0\\
0.125&0,70710& 0,70485&0,70629 &0,70681&0,70700\\
0.250&1,0& 0,99793& 0,99927&0,99974 &0,99991 \\
0.50&1,41421&1,41244 &1,41359&1,41399&1,41413\\
0.75& 1,73205& 1,73050& 1,73151&1,73186 &1,73198\\
1.0&2.0 &1,99862 &1,99952&1,99983&1,99994\\ \hline
\end{tabular}
\vspace{0.3cm}

Table 1. The approximative solution for $N=16, 32, 64, 128.$
\vspace{0.3cm}

\noindent In this table: $y(x)$-- the sharp solution, $\tilde y(x)$--the approximative solution.\\
{\bf Example 4.} Let us consider the Cauchy problem

\begin{equation} \label{eq5.2}
\left\{
\begin{array}{ll}
y'+p'(x)y=q'(x),&\quad x\in [0,1]\\
y(0)=0.
\end{array}
\right.
\end{equation}
where
$$
p(x)=\left\{
      \begin{array}{lcl}
      \sqrt{x}&if&x\in[0,1/3],\\
      (2/3-x)\sqrt{3}&if&x\in[1/3,2/3],\\
      \sqrt{x-2/3}&if&x\in[2/3,1],\\
      \end{array}
      \right.
$$
$$
p'(x)=\left\{
      \begin{array}{lcl}
      \frac{1}{2\sqrt{x}}&if&x\in[0,1/3],\\
      -\sqrt{3}&if&x\in[1/3,2/3],\\
      \frac{1}{2\sqrt{x-2/3}}&if&x\in[2/3,1],\\
      \end{array}
      \right.
$$
$$
q(x)=\left\{
      \begin{array}{lcl}
      \frac{x}{2}+ \sqrt{x}-2/3&if&x\in[0,1/3],\\
      -x(2+\sqrt{3})+\frac32 x^2+\frac{2}{\sqrt{3}} &if&x\in[1/3,2/3],\\
      \frac{x}{2} +\sqrt{x-2/3}-1&xif&\in[2/3,1],\\
      \end{array}
      \right.
$$
$$
q'(x)=\left\{
      \begin{array}{lcl}
      \frac12 (1 +\frac{1}{\sqrt{x}})&if&x\in[0,1/3],\\
      -2-\sqrt{3}+3x&if&x\in[1/3,2/3].\\
      \frac{1}{2}(1 +\frac{1}{\sqrt{x-2/3}})&if&x\in[2/3,1].\\
      \end{array}
      \right.
$$
The solution
$$
y(x)=\left\{
      \begin{array}{ll}
      \sqrt{x}&x\in[0,1/3]\\
      \sqrt{3}(2/3-x)&x\in[1/3,2/3]\\
      \sqrt{x-2/3}&x\in[2/3,1].\\
      \end{array}
      \right.
$$
  of this problem is continuous function on $[0,1]$, but the derivatives $y'(\frac13),\ y'(\frac23),\ y'(0)$ not exist.
  In Figure 1 we give graphs of the approximate (blue) and exact (red) solutions. Both graphs are drown  on 512 points. 
We see that these graphs are the same.

\vspace{0.2cm}

\hspace{0.5cm}\includegraphics[scale=0.33]{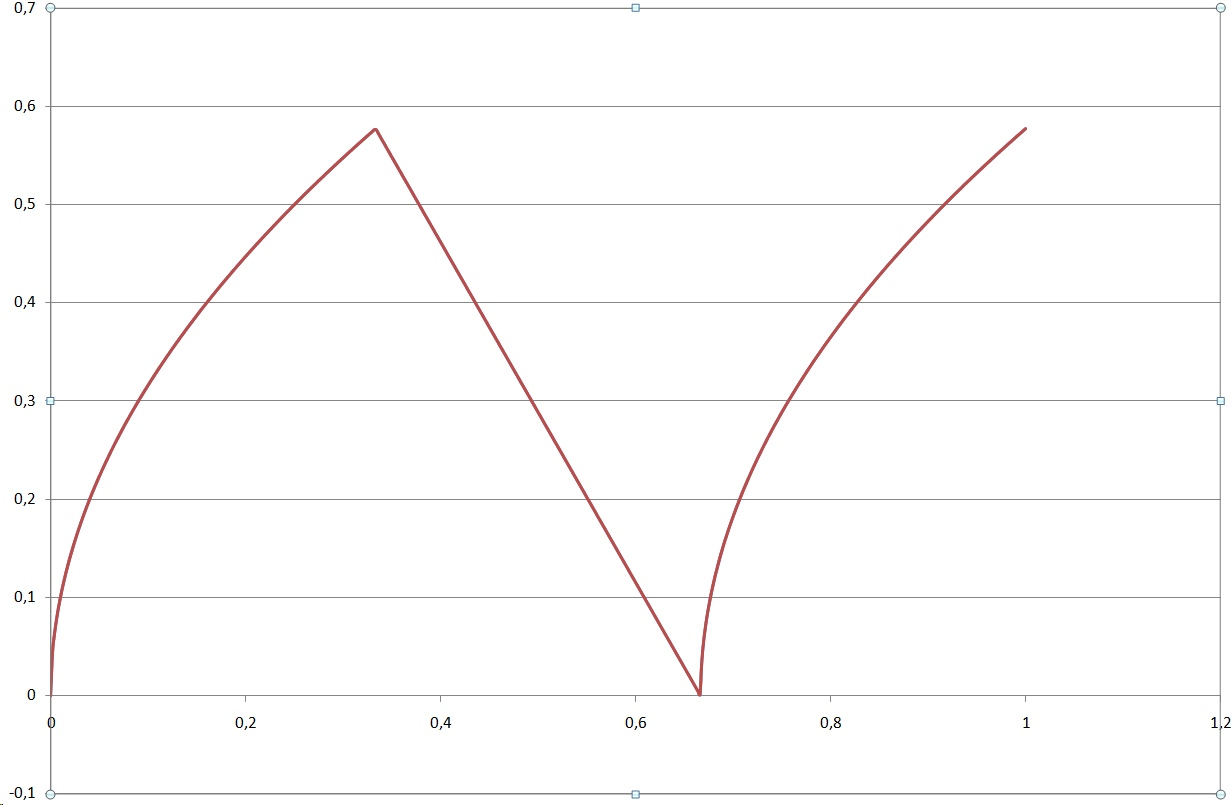}

\vspace{0.3cm}
\hspace{0.8cm}Figure 1. The graphs of  $\tilde y(x)$-(blue) and $y(x)$-(red )for $2^n=N=512$.

\vspace{0.3cm}
\noindent
   Denote by  $\tilde y_n(x)$ the approximative solution for the point system $(j2^{-n})_{j=0}^{2^n}$ and
 $\delta_n=\max_j |\tilde y_n(j2^{-n})-y(j2^{-n})| $.
 In the Table 2 we give the error of the approximative solution of Cauchy problem (\ref{eq5.2}) for $n=\overline{4,10}$.
\vspace{0.3cm}

\hspace{0.5cm}\begin{tabular}{|l|r|r|r|r|r|r|r|}\hline
&n=4&n=5&n=6&n=7&n=8&n=9&n=10\\ \hline
$\delta_n$&$1.1\cdot10^{-3}$&$5.3\cdot10^{-4}$&$1.8\cdot10^{-4}$
&$8.6\cdot10^{-5}$&$2.8\cdot10^{-5}$&$1.2\cdot10^{-5}$&$3.9\cdot10^{-6}$\\ \hline
\end{tabular}
\vspace{0.3cm}

\hspace{0.5cm}Table 2. The error of the approximative solution for $2^n=16, 32, 64, 128,256,512,1024.$
\vspace{0.3cm}

Acknowledgements.This work was supported by SEMC "Mathematics of Future Technologies".

 \end{document}